\documentclass[a4paper,12pt]{article}
\usepackage[utf8]{inputenc}
\usepackage[T1]{fontenc}
\usepackage{amsmath}
\usepackage{amsthm}
\usepackage{amssymb}
\usepackage{color}
\usepackage[all]{xy}
\usepackage{setspace}
\usepackage{mathtools}
\usepackage{enumitem}

\title{Global algebras of\\ nonlinear generalized functions with\\ applications in general relativity}
\author{Eduard A.~Nigsch\footnote{Fakultät für Mathematik, Universität Wien, Oskar-Morgenstern-Platz 1, 1090 Wien, Austria; e-mail: eduard.nigsch@univie.ac.at; phone: +43 1 4277 50760}, Clemens Sämann\footnote{Fakultät für Mathematik, Universität Wien, Oskar-Morgenstern-Platz 1, 1090 Wien, Austria; e-mail: clemens.saemann@univie.ac.at; phone: +43 1 4277 50632}}

\newcommand{\D}{\mathrm{D}}
\newcommand{\ud}{\mathrm{d}}
\newcommand{\cD}{\mathcal{D}}
\newcommand{\bR}{\mathbb{R}}
\newcommand{\bC}{\mathbb{C}}
\newcommand{\bN}{\mathbb{N}}
\newcommand{\cA}{\mathcal{A}}
\newcommand{\cE}{\mathcal{E}}
\newcommand{\cG}{\mathcal{G}}
\newcommand{\cS}{\mathcal{S}}
\newcommand{\cT}{\mathcal{T}}
\newcommand{\cL}{\mathcal{L}}

\newcommand{\cN}{\mathcal{N}}
\newcommand{\pder}[2]{\frac{\pd#1}{\pd#2}}
\newcommand{\vp}{\mathrm{vp}}
\newcommand{\fX}{\mathfrak{X}}
\newcommand{\coleq}{\coloneqq}
\newcommand{\e}{\varepsilon}
\newcommand{\abso}[1]{\lvert#1\rvert}
\newcommand{\norm}[1]{\lVert#1\rVert}

\newcommand{\pd}{\partial}
\newcommand{\csub}{\subset\subset}
\DeclareMathOperator{\id}{id}
\DeclareMathOperator{\supp}{supp}
\DeclareMathOperator{\cs}{cs}
\newcommand{\Lie}{\mathrm{L}}

\newcommand{\Gob}[1]{\cE^o(#1)}
\newcommand{\Gom}[1]{\cE^o_M(#1)}
\newcommand{\Gon}[1]{\cN^o(#1)}
\newcommand{\Go}[1]{\cG^o(#1)}
\newcommand{\SGo}{\cG^o}
\newcommand{\Geb}[1]{\cE^e(#1)}
\newcommand{\Gem}[1]{\cE^e_M(#1)}
\newcommand{\Gen}[1]{\cN^e(#1)}
\newcommand{\Ge}[1]{\cG^e(#1)}
\newcommand{\SGe}{\cG^e}
\newcommand{\Gsb}[1]{\cE^s(#1)}
\newcommand{\Gsm}[1]{\cE^s_M(#1)}
\newcommand{\Gsn}[1]{\cN^s(#1)}
\newcommand{\Gs}[1]{\cG^s(#1)}
\newcommand{\SGs}{\cG^s}
\newcommand{\Gdb}[1]{\cE^d(#1)}
\newcommand{\Gdm}[1]{\cE^d_M(#1)}
\newcommand{\Gdn}[1]{\cN^d(#1)}
\newcommand{\Gd}[1]{\cG^d(#1)}
\newcommand{\SGd}{\cG^d}
\newcommand{\Gnb}[1]{\cE(#1)}
\newcommand{\Gnm}[1]{\cE_M(#1)}
\newcommand{\Gnn}[1]{\cN(#1)}
\newcommand{\Gn}[1]{\cG(#1)}
\newcommand{\SGn}{\cG}
\newcommand{\Ghb}[1]{\hat{\cE}(#1)}
\newcommand{\Ghm}[1]{\hat{\cE}_M(#1)}
\newcommand{\Ghn}[1]{\hat{\cN}(#1)}
\newcommand{\Gh}[1]{\hat{\cG}(#1)}
\newcommand{\SGh}{\hat{\cG}}
\newcommand{\Ghrs}[1]{\hat\cG^r_s(#1)}
\newcommand{\SGhrs}{\hat{\cG}^r_s}
\newcommand{\Gfb}[1]{\cE(#1)}
\newcommand{\Gflb}[1]{\cL(#1)}
\newcommand{\Gfbloc}[1]{\cE_{\mathrm{loc}}(#1)}
\newcommand{\SGfbloc}{\cE_{\mathrm{loc}}}
\newcommand{\Gfbploc}[1]{\cE_{\mathrm{ploc}}(#1)}
\newcommand{\SGfbploc}{\cE_{\mathrm{ploc}}}
\newcommand{\Gfbpi}[1]{\cE_{\mathrm{pi}}(#1)}
\newcommand{\SGfbpi}{\cE_{\mathrm{pi}}}
\newcommand{\Gflbloc}[1]{\cL_{\mathrm{loc}}(#1)}

\newcommand{\Gflbploc}[1]{\cL_{\mathrm{ploc}}(#1)}
\newcommand{\Gflbpi}[1]{\cL_{\mathrm{pi}}(#1)}

\newcommand{\SK}[1]{\mathrm{SK}(#1)}

\newtheorem{theorem}{Theorem}
\newtheorem{definition}[theorem]{Definition}
\newtheorem{proposition}[theorem]{Proposition}

\begin{document}

\maketitle

\begin{abstract}
We give an overview of the development of algebras of generalized functions in the sense of Colombeau and recent advances concerning diffeomorphism invariant global algebras of generalized functions and tensor fields. We furthermore provide a survey on possible applications in general relativity in light of the limitations of distribution theory.
\end{abstract}



\section{Introduction}\label{sec_intro}

In physics there is the natural need to reduce the complexity of real-world physical systems for the formulation of mathematical models and therefore to work with idealizations. A well-known example of this is the notion of a point particle, which represents an object whose properties like mass or electric charge are concentrated in a single point.

Such a point charge is classically represented by Dirac's delta function $\delta$, originally introduced as a useful notation for dealing with a certain kind of infinities \cite[p.~58]{zbMATH03131829}:
\begin{quote}
 ``To get a precise notation for dealing with these infinities, we introduce a quantity $\delta(x)$ depending on a parameter $x$ satisfying the conditions
$\int_{-\infty}^\infty \delta(x)\,\ud x = 1$, $\delta(x) = 0\ \textrm{ for }x \ne 0$. [...] $\delta(x)$ is not a function according to the usual mathematical definition of a function, which requires a function to have a definite value for each point in its domain, but is something more general, which we may call an `improper function' to show up its difference from a function defined by the usual definition. Thus $\delta(x)$ is not a quantity which can be generally used in mathematical analysis like an ordinary function, but its use must be confined to certain simple types of expression for which it is obvious that no inconsistency can arise.``
\end{quote}

Its mathematical justification had to wait for the introduction of the theory of distributions by S.~L.~Sobolev \cite{0014.05902} and L.~Schwartz \cite{MR0035918}, who defined distributions as continuous linear functionals on certain spaces of test functions.

At that time the theory of distributions not only furnished the means for a rigorous formulation of several previously vague uses of generalized functions in physics, but very quickly brought strong results in particular in the field of differential equations. As the most prominent example, the theorem of Malgrange-Ehrenpreis \cite{zbMATH03118689,zbMATH03089436,wagnermalgrange} states that every linear partial differential operator with constant coefficients has a distributional fundamental solution. Another result of foundational importance is Schwartz' kernel theorem \cite{FDVV}, which allows one to represent operators between general spaces of distributions by distributional kernels. These results reflect and in effect are made possible by the strong footing of distribution theory in the theory of locally convex spaces. At the same time, it is noteworthy that the basic definitions of distribution theory are very simple, which certainly makes it more accessible also for non-specialists and partly explains its success in many fields.

It is clear from the quote above that there are certain restrictions on what one can do with distributions compared with ordinary functions. Being an inherently linear theory, distribution theory does not allow for an intrinsic definition of an associative product of distributions preserving the pointwise product of continuous functions, as can be seen from the relations
\[ 0=(\delta(x)\cdot x)\cdot \vp\left(\frac{1}{x}\right) = \delta(x)\cdot (x\cdot \vp\left(\frac{1}{x}\right)) =
\delta(x),\]
since $\delta(x)\cdot x=x|_{x=0}=0$ and $x\cdot\vp({1}/{x}) = 1$, where $\vp(1/x)$ denotes the \emph{Cauchy principal
value of ${1}/{x}$}. A theory where $\delta=0$ would be pointless, since we explicitly want to model shocks and point
charges. Similarly, one can see that commutativity and the Leibniz rules cannot hold simultaneously. 
This impossibility is stated in precise terms by the Schwartz impossibility result \cite{Schwartz}:
\begin{theorem}\label{schwartz}Let $\cA$ be an algebra containing the algebra $C(\bR)$ of all continuous functions on $\bR$ as a subalgebra such that the constant function $1 \in C(\bR)$ is the unit in $\cA$. Assume there exists a linear map $D \colon \cA \to \cA$ extending the derivation of continuously differentiable functions and satisfying the Leibnitz rule. Then $D^2 ( \abso{x} ) = 0$.
\end{theorem}
Because $D^2(\abso{x}) = 2\delta$ this implies that $\cA$ cannot contain $\delta$ (for a further discussion of the problem of multiplying distributions see \cite[Chapter I]{MOBook}). Hence, one way to define a product of distributions is to restrict to subalgebras (e.g., Sobolev spaces $H^s$ with $s$ large enough) or certain pairs of distributions (e.g., $L^p \times L^q$ with $1/p+1/q=1$) for which the product again makes sense as a distribution; a different way is to weaken the requirement that $C(\bR)$ is a subalgebra.

Both approaches have their merits. In Section \ref{sec_lindist} we first look at the uses and limitations of the first approach in the context of distributional geometry, where one only considers distributions such that all desired calculations can be carried out. After that, we will examine algebras of nonlinear generalized functions in the sense of Colombeau, which contain the space of smooth functions as a subalgebra. After recalling the basic idea and original definitions in Section \ref{sec_colalg} we describe how they split up into two different variants in Section \ref{sec_fullspec} and present some recent applications to general relativity in Section \ref{sec_appl}. In Section \ref{sec_diffinv} we then describe the previous attempts on obtaining diffeomorphism invariant Colombeau algebras, including the tensor case. Finally, in Section \ref{sec_funcanal} we present a new approach to Colombeau algebras, unifying and simplifying previous work and leading the way to further developments.

\section{Linear distributional geometry}\label{sec_lindist}

\subsection{Basic definitions}

We briefly recall the basic definitions of the theory of distributions on an open subset $\Omega \subseteq \bR^n$. For introductory texts we refer to \cite{Friedlander,MR0035918,zbMATH05622635,zbMATH03230708,Treves}.

The space of \emph{test functions} $\cD(\Omega)$ consists of all functions $\varphi\colon \Omega \to \bC$ which are smooth (infinitely differentiable) and have compact support. This space is endowed with a certain locally convex inductive limit topology which turns it into an LF-space. The space of distributions $\cD'(\Omega)$ then is defined as its topological dual endowed with the strong topology. The action of a distribution $u$ on a test function $\varphi$ is commonly denoted by $\langle u, \varphi \rangle$.

Locally integrable functions $f$ can be embedded into $\cD'(\Omega)$ via integration; the image of $f$ in $\cD'(\Omega)$ then is called a \emph{regular} distribution and is given by the functional
\begin{equation}\label{regdist}
\varphi \mapsto \int f(x) \varphi(x)\,\ud x\qquad (\varphi \in \cD(\Omega)).
\end{equation}
This formula also provides the basis on which many classical operations are extended to distributions: replacing $f$ by its $i$-th partial derivative $\pd_i f$ and integrating by parts, or replacing $f$ by its product $g \cdot f$ with a smooth function $g$, one obtains from \eqref{regdist} formulas for differentiating a distribution $u$ or multiplying it by a smooth function:
\begin{align*}
\langle \pd_i u, \varphi \rangle & \coleq - \langle u, \pd_i \varphi \rangle, \\
\langle g \cdot u, \varphi \rangle & \coleq \langle u, g \cdot \varphi \rangle.
\end{align*}

If one wants to extend the notion of distributions to a manifold $M$ (which will always be assumed to be finite-dimensional, paracompact and orientable), the definition of regular distributions as in \eqref{regdist} requires the product of $f$ and $\varphi$ to be an $n$-form, as these are the objects which can be integrated on a manifold. For this reason, one takes \emph{compactly supported $n$-forms} $\omega \in \Omega^n_c(M)$ instead of test functions and defines $\cD'(M)$ as the dual of $\Omega^n_c(M)$. As above, the Lie derivative of functions extends to $\cD'(M)$ by setting $\langle \Lie_X u, \omega \rangle \coleq - \langle u, \Lie_X \omega \rangle$ for $u \in \cD'(M)$, $\omega \in \Omega^n_c(M)$ and $X$ a smooth vector field on $M$.

For considering e.g.~distributional sources in linear field theories or singular metrics in general relativity one needs a notion of distributional sections of vector bundles or, more specifically, distributional tensor fields. The latter are most easily defined as tensor fields with distributional coefficients, hence we say that distributional tensors of rank $(r,s)$ are given by the tensor product
\[ \cD'^r_s(M) \coleq \cD'(M) \otimes_{C^\infty(M)} \cT^r_s(M), \]
where $\cT^r_s(M)$ is the space of smooth $(r,s)$-tensor fields on $M$ and $C^\infty(M)$ the space of smooth functions on $M$ (see \cite{GKOS} for a thorough introduction to vector bundle valued distributions on manifolds).

Trying to extend classical operations to distributional tensor fields one encounters already the first difficulty: applying the method of continuous extension to multilinear operations on smooth tensor fields allows at most one factor to be distributional, the others have to be smooth. For instance, one can only take the tensor product of a distributional tensor field with a smooth tensor field, but not with another distributional tensor field. Furthermore, this places restrictions on the definition of distributional metrics and connections (\cite[Section 3.1.5]{GKOS}). The underlying reason can be seen from the respective coordinate expressions, which would involve multiplication of distributions.

\subsection{The Geroch-Traschen class}

Geroch and Traschen \cite{gerochtraschen} tried to find a wide class of metrics whose curvature tensors make sense as distributions --- the so-called \emph{Geroch-Traschen class of metrics} (or \emph{gt-regular metrics} in \cite{SV:09}). For this purpose they defined certain restrictions on the metric and the covariant derivative in order to make all desired quantities well-defined.

The coordinate expression of the curvature tensor $R$ of a semi-Riemannian manifold $(M,g)$ reads
\begin{equation*}
 R^i_{jkl} = \pder{\Gamma^i_{kj}}{x^l} - \pder{\Gamma^i_{lj}}{x^k} + \sum_{m=1}^n \Gamma^i_{lm} \Gamma^m_{kj} - \sum_{m=1}^n
\Gamma^i_{km}\Gamma^m_{lj},
\end{equation*}
where the Christoffel symbols $\Gamma^k_{ij}$ are given by
\[ \Gamma^k_{ij}=\frac{1}{2}\sum_{m=1}^n g^{km}\bigl( \pder{g_{jm}}{x^i} +
\pder{g_{im}}{x^j} - \pder{g_{ij}}{x^m} \bigr). \]
Here, $g_{ij}$ denotes the components of the metric and $g^{ij}$ those of its inverse. The crucial observation of Geroch and Traschen was that the following conditions are sufficient for $R^i_{jkl}$ to define a distribution:
\begin{enumerate}[label=(\roman*)]
 \item The inverse of $g$ exists everywhere and $g_{ij}, g^{ij}$ are locally bounded.
 \item The weak first derivative of $g_{ij}$ exists and is locally square integrable.
\end{enumerate}
This conditions are however not necessary, as can be seen in \cite{G:99}. On the other hand the broader class of metrics
given in \cite{G:99} (yielding distributional curvature) is not stable under appropriate approximations via smooth metrics.
More on smoothing of gt-regular metrics, non-degeneracy and compatibility with the Colombeau approach can be found in
\cite{SV:09}. It should be noted that there is also a coordinate-free way to derive the Geroch-Traschen class of metrics
(\cite{leFMlorentzian}).

\subsection{Limitations}\label{sec_limits}

The most important result in the context of applicability of distribution theory to general relativity is \cite[Theorem 1]{gerochtraschen}: if one has a gt-regular distribution $u\neq 0$ concentrated on a submanifold $S$ (i.e.,
$\text{supp}(u)\subseteq S$) then the codimension of $S$ has to be equal to one. Thus in four dimensions this excludes strings and points
and only allows concentration of matter on shells, which suggests that the use of distributions in general relativity is quite
limited. We will discuss this in some more detail. 

The application of distribution theory in general relativity has a long history, for example with problems involving metrics of low regularity, including shocks and thin shells (cf.~\cite{CD:87,Israel:66,Israel:67,KP:71, B:97, AB:96, BN:94} and \cite{genrel} for a survey). As emphasized in \cite{genrel} this approach is severely limited in its range of applications since on one hand the theory of distributions is a linear theory and on the other hand general relativity is inherently nonlinear (due to Einstein's field equations). We will exemplify this problem in the case of impulsive pp-waves (plane fronted waves with parallel rays), where the nonlinearity comes from the geodesic equations.

\subsection{Impulsive pp-waves}
Impulsive waves were introduced by Penrose (\cite{Penrose:68, Penrose:72}) using a ``cut-and-paste'' approach, which involves cutting Minkowski space-time along a null hypersurface and reattaching the two regions with a suitable warp. These waves can also be considered as idealizations (impulsive limit) of sandwich waves of infinitely short duration. We refer to \cite[Chapter 20]{GrPo:09} for a general introduction to this topic. Impulsive pp-waves can be described by the so-called \emph{Brinkmann
form} of the pp-wave metric (\cite{Brinkmann:23})
\begin{equation}\label{eq-ipp}
 \ud s^2 = f(x,y)\delta(u) \ud u^2 - \ud u \ud v + \ud x^2 + \ud y^2,
\end{equation}
where $f\in C^\infty(\mathbb{R}^2)$, $\delta$ denotes the delta distribution, $u,v$ are null coordinates such that $\pd_v$ is covariantly constant and $x,y$ are Cartesian coordinates spanning the wave surfaces. From this form of the metric one can easily see that the space-time is flat except on the null hyperplane $u=0$, where a $\delta$-like impulse is located. To be more precise, the impulse has ``strength'' $f(x,y)$ at the point $(x,y,0,v)$. Clearly this space-time cannot be described within classical general relativity due to the $\delta$-like impulse. Moreover, it is not even gt-regular, i.e., it does not belong to the Geroch-Traschen class of metrics. Nevertheless these space-times are ``tame'' enough to calculate the Ricci tensor within distribution theory and hence the Einstein vacuum equations may be formulated giving $\Delta f=0$ on the hypersurface $u=0$.

So what aspects of this space-time cannot be handled in distribution theory? In light of the problem of multiplying distributions we see that when trying to solve the geodesic equations for
impulsive pp-wave space-times described by the metric \eqref{eq-ipp}, we are in general multiplying distributions, not
classical functions. To describe this problem in more detail, we consider the geodesic equations for the
metric \eqref{eq-ipp}. Since $\ddot{u}=0$, we use $u$ as a new affine parameter and thereby exclude only trivial geodesics
parallel to the shock hypersurface. In summary we get
\begin{align}
\ddot{v}(u)&=f(u) \dot{\delta}(u) + 2\left(\pder{f}{x}(u) x(u) + \pder{f}{y}(u) y(u)\right)\delta(u)\label{geo-v},\\
\ddot{x}(u)&=\frac{1}{2}\pder{f}{x}(u) \delta(u),\\
\ddot{y}(u)&=\frac{1}{2}\pder{f}{y}(u) \delta(u).
\end{align}
We try to solve these equations in $\cD'$: for $u\neq 0$ the right hand sides vanish, so we expect that the
geodesics are broken (and possibly refracted) straight lines. By integrating and simplifying we get an expression for $x$, namely $x(u)=x_0 + \frac{1}{2}(\pd f / \pd x)(x_0)H(u)u$, where $H$ denotes the Heaviside function and $x_0=x(0)$ is an arbitrary
initial condition. Consequently, in equation \eqref{geo-v} we have the product $H\cdot\delta$ which is
ill-defined in classical distribution theory. Now we see that the problem is that the distributional ``solutions'' do not
obey the original equations unless one is willing to impose certain ad hoc multiplication rules. The correct idea to solve
this problem is to approximate the delta distribution by a net of smooth functions $(\delta_\e)_\e$ ($\e \in (0,1]$)
and solve the equations in $C^\infty$ for each $\e$. This idea can be made mathematically rigorous in the Colombeau setting.
See \cite{Steinbauer:98, geodesics,KS:99} for a treatment of impulsive pp-waves in the context of Colombeau algebras.

\section{Algebras of nonlinear generalized functions}\label{sec_colalg}

Theorem \ref{schwartz} implies that one cannot hope to embed distributions into a differential algebra such that the product of continuous functions is preserved. However, it was discovered by J.~F.~Colombeau \cite{ColElem, ColNew} that one can construct associative commutative differential algebras containing $\cD'$ as a linear subspace and $C^\infty$ as a subalgebra. In the following we will briefly sketch his original ideas.

Colombeau's algebras emerged in the context of his work on calculus in infinite dimensional locally convex spaces. Clearly, products of distributions on an open subset $\Omega \subseteq \bR^n$ can be naturally seen as monomials on $\cD(\Omega)$ (in the sense of \cite{dineen}) or, more generally, smooth functions on $\cD(\Omega)$. The notion of smoothness employed by Colombeau originally was that of Silva-differentiability (\cite{zbMATH03798426}) but can be replaced with the (here equivalent) notion of smoothness in the sense of calculus on convenient vector spaces (\cite{KM}), which we shall also use in the sequel because it appears to be somewhat simpler for our purposes.

One might define the product of $u,v \in \cD'(\Omega)$ as the mapping $\cD(\Omega) \to \bC$ given by
\[ (u \cdot v)(\varphi) \coleq \langle u, \varphi \rangle \cdot \langle v, \varphi \rangle \qquad (\varphi \in \cD(\Omega)). \] 
Obviously $C^\infty(\cD(\Omega))$, the space of smooth functions from $\cD(\Omega)$ into $\bC$, is an algebra containing $\cD'(\Omega)$. It possesses partial derivatives $\pd_i$ ($i=1 \dotsc n$) extending those of distributions by setting, for $R \in C^\infty(\cD(\Omega))$,
\begin{equation}\label{star}
(\pd_i R)(\varphi) \coleq - \ud R (\varphi) \cdot \pd_i \varphi \qquad (\varphi \in \cD(\Omega))
\end{equation}
where $\ud R$ denotes the differential of $R$ (\cite[3.18]{KM}). The product of smooth functions $f,g$ is not preserved in this algebra because the expression
\begin{equation}\label{blah}
\int f(x) \varphi(x) \,\ud x \cdot \int g(x) \varphi(x) \,\ud x - \int f(x)g(x) \varphi(x) \,\ud x
\end{equation}
is nonzero in general. It vanishes, however, if $\varphi$ is replaced by $\delta_z$, the delta distribution at any point $z \in \Omega$. Colombeau's brilliant idea here was to characterize elements of the form \eqref{blah} in $C^\infty(\cD(\Omega))$ by their asymptotic behavior on certain sequences of test functions and hence define an ideal containing them. The corresponding quotient algebra would then preserve the product of smooth functions. We will explain this idea in some more detail because of its fundamental importance for the rest of this article; for the original account, see \cite{ColNew}.

First, $C^\infty(\Omega)$ as a reflexive space is topologically isomorphic to its bidual $L_b( \cE'(\Omega), \bC)$, where $\cE'(\Omega)$ denotes the strong dual of $C^\infty(\Omega)$ and, given locally convex spaces $E$ and $F$, $L_b(E,F)$ denotes the space of all continuous linear mappings from $E$ to $F$ endowed with the topology of bounded convergence. We see that the subalgebra of $C^\infty(\cD(\Omega))$ generated by $C^\infty(\Omega)$ is contained in $C^\infty( \cE'(\Omega))$:
\[
 \xymatrix{
\cD'(\Omega) \ar[rr] && C^\infty(\cD(\Omega)) \\
C^\infty(\Omega) \ar[u] \ar@{}[r]|-*[@]\txt{$\cong$} & L(\cE'(\Omega), \bC) \ar@{}[r]|-*[@]\txt{$\subseteq$} & C^\infty(\cE'(\Omega)) \ar[u]
}
\]
Slightly generalizing \eqref{blah}, we can say that we are looking for an ideal of $C^\infty(\cD(\Omega))$ containing the set
\begin{equation}\label{alpha}
\{ R \in C^\infty ( \cE'(\Omega) )\ |\ R ( \delta_x ) = 0\ \forall x \in \Omega \}.
\end{equation}
The following proposition characterizes elements of this set in terms of their restriction to $\cD(\Omega)$ (\cite[Proposition 3.3.3]{ColNew}). It utilizes the sets
\begin{multline}\label{locmoment}
 \cA_q(\bR^n) \coleq \{ \varphi \in \cD(\bR^n): \int \varphi(x)\,\ud x = 1,\\
 \int x^\alpha \varphi(x)\,\ud x = 0\ \forall\ 1 \le \abso{\alpha} \le q \} \qquad (q \in \bN_0)
\end{multline}
using the usual multiindex notation $x^\alpha = x_1^{\alpha_1}\dotsm x_n^{\alpha_n}$ for $x=(x_1,\dotsc,x_n) \in \bR^n$ and $\alpha = (\alpha_1,\dotsc,\alpha_n) \in \bN_0^n$.
\begin{proposition}
 Let $R \in C^\infty ( \cE'(\Omega) )$ be given. Then $R ( \delta_x ) = 0$ for all $x \in \Omega$ if and only if $\forall q \in \bN$ $\forall \varphi \in \cA_q(\bR^n)$ $\forall K \subseteq \Omega$ compact there are $c>0$ and $\eta > 0$ such that
\[ \abso{ R(\varphi_{\e,x}) } \le c \e^{q+1}\quad \forall\ 0 < \e < \eta\textrm{ and }x \in K, \]
where $\varphi_{\e,x} \in \cD(\Omega)$ is defined by $\varphi_{\e,x}(y) \coleq \e^{-n} \varphi( (y-x) / \e )$ for small $\e$.
\end{proposition}
This suggests to define an ideal of $C^\infty(\cD(\Omega))$ containing \eqref{alpha} as the set of all elements satisfying this condition. However, one has to restrict to a smaller subalgebra (still containing all distributions) for this to be an ideal. The definition of the original algebra of nonlinear generalized functions introduced by Colombeau hence takes the following form.

\begin{definition}\label{deftest}If $R \in C^\infty(\cD(\Omega))$ we say that $R$ is moderate if for every compact subset $K$ of $\Omega$ and every multiindex $\alpha \in \bN_0^n$ there is an $N \in \bN$ such that $\forall \varphi \in \cA_N(\bR^n)$ $\exists c>0$ and $\eta >0$ such that $\abso{\pd^\alpha R(\varphi_{\e,x})} \le c \e^{-N}$. The set of all moderate elements is denoted by $\Gom\Omega$.
 
$R \in \Gom\Omega$ is called negligible if for every compact set $K$ of $\Omega$ and every multiindex $\alpha \in \bN_0^n$ there is an $N \in \bN$ such that $\forall \varphi \in \cA_q(\bR^n)$ with $q \ge N$ $\exists c>0$ and $\eta>0$ such that $\abso{\pd^\alpha R(\varphi_{\e,x})} \le c \e^{q-N}$. The subset of all negligible elements is denoted by $\Gon\Omega$.

We set $\Gob\Omega \coleq C^\infty(\cD(\Omega))$ and $\Go\Omega \coleq \Gom\Omega / \Gon\Omega$. 
\end{definition}

For clarity, let us state that $\pd^\alpha R(\varphi_{\e,x})$ denotes the derivative of order $\alpha$ of the map $x' \mapsto R(\varphi_{\e,x'})$.

The desired properties are easily verified: distributions are canonically embedded into $\Go\Omega$, $C^\infty(\Omega)$ is contained as a subalgebra and the partial derivatives given by \eqref{star} extend those of distributions. Furthermore, many concepts of classical analysis (point values, integrals, Fourier transform etc.) have direct equivalents in the context of generalized functions (cf.~\cite{ColNew,ColElem,MOBook,GKOS}).

It should be stressed at this point that Colombeau's algebra of generalized functions is not an alternative to but an extension of distribution theory. In fact, through the concept of \emph{association} it is possible in many cases to relate certain generalized functions in $\Go\Omega$ to distributions (\cite[Definition 3.5.2]{ColNew}):
\begin{definition}Let $R \in \Gob\Omega$ be a representative of an element of $\Go\Omega$. If for every $\psi \in \cD(\Omega)$ there is $q \in \bN$ such that for all $\varphi \in \cA_q(\bR^n)$ the limit
\[ \lim_{\e \to 0} \int R(\varphi_{\e,x}) \cdot \psi(x)\,\ud x \]
exists independently of $\varphi$, and as a function of $\psi$ this limit defines a distribution on $\Omega$, we say that $R$ admits an \emph{associated distribution} $u$ defined by the formula
\[ \langle u, \psi \rangle \coleq \lim_{\e \to 0} \int R(\varphi_{\e,x}) \cdot \psi(x)\,\ud x\qquad (\psi \in \cD(\Omega)). \]
\end{definition}

As an example, the classical product of continuous functions is not preserved in $\Go\Omega$, but it is so on the level of association. The point of association is that whenever a calculation makes sense distributionally one obtains the same result if one formulates the same problem in the Colombeau algebra, does the same calculations there, and finds the associated distribution.

These definitions are classical and have since been subject to many modifications and adaptions, of which we will describe the most important ones in Section \ref{sec_fullspec}. For the many applications of Colombeau algebras in the fields of partial differential equations, numerics or geometry we refer to the monographs \cite{ColNew, Biagioni, MOBook, GKOS}.

\section{Full and special Colombeau algebras}\label{sec_fullspec}

After the presentation of $\SGo$ in \cite{ColNew}, two main variants of Colombeau algebras emerged. The first one, so-called \emph{full} algebras, possess a canonical embedding of distributions and were the basis for the development of a diffeomorphism invariant theory which eventually should lead to a global formulation on manifolds and also an algebra of generalized tensor fields. The second one, so-called \emph{simplified} or \emph{special} algebras, have a very simple formulation and thus allow for a direct transfer of many classical notions, but have no canonical embedding of distributions anymore and suffer from a lack of diffeomorphism invariance.

At the basis of these variants lie different interpretations of generalized functions in $\Go\Omega$ on a rather symbolic level. In fact, by Definition \ref{deftest} generalized functions $R \in \Gob\Omega$ are determined by all their values
\begin{equation}\label{genfunct}
R(\varphi_{\e,x})\qquad (\varphi \in \cA_1(\bR^n),\ \e \in (0,1],\ x \in \Omega).
\end{equation}
By varying how the dependence on $\varphi$, $\e$ and $x$ is incorporated in the basic space the following prototypical variants can be derived.

\subsection{Full algebras}

The idea of the \emph{elementary full} algebra presented in \cite{ColElem} is to view generalized functions as mappings in both $\varphi$ and $x$. This allows one to shift the burden of differentiation to the $x$-variable alone, for fixed $\varphi$, and formulate the theory without recourse to calculus on infinite dimensional locally convex spaces. The corresponding definitions are as follows (see \cite[Chapter 1]{ColElem} and \cite[Section 1.4]{GKOS}).
\begin{definition}We set
\begin{align*}
 U(\Omega) &\coleq \{(\varphi,x) \in \cA_0(\bR^n) \times \Omega\ |\ \supp \varphi +x \subseteq \Omega \}, \\
\Geb{\Omega} &\coleq \{ R \colon U(\Omega) \to \bC\ |\ R(\varphi,.)\textrm{ is smooth }\forall \varphi \in \cA_0(\bR^n)\}.
\end{align*}
We say that an element $R \in \Geb\Omega$ is moderate if for every compact subset $K$ of $\Omega$ and every multiindex $\alpha \in \bN_0^n$ there is $N \in \bN$ such that if $\varphi \in \cA_N(\bR^n)$ then there exist $\eta>0$, $c>0$ such that $\abso{\pd^\alpha R(S_\e \varphi, x)} \le c \e^{-N}$, where $(S_\e \varphi)(y) \coleq \e^{-n} \varphi(y/\e)$. We denote by $\Gem\Omega$ the set of all moderate elements of $\Geb\Omega$.

We say that an element $R$ of $\Geb\Omega$ is negligible if for every compact subset $K$ of $\Omega$, every multiindex $\alpha \in \bN_0^n$ and every $m \in \bN$ there is $q \in \bN$ such that if $\varphi \in \cA_q(\bR^n)$ then there exist $\eta>0, c>0$ such that $\abso{\pd^\alpha R(\varphi_\e, x)} \le c \e^q$. We denote by $\Gen\Omega$ the set of all negligible elements.

We set $\Ge\Omega \coleq \Gem\Omega / \Gen\Omega$.
\end{definition}

The canonical embedding $\iota\colon\cD'(\Omega) \to \Ge\Omega$ is determined by assigning to $u \in \cD'(\Omega)$ the class of
\[ (\varphi, x) \mapsto \langle u, \varphi(.-x) \rangle\qquad ((\varphi,x) \in U(\Omega)). \]
Smooth functions are embedded via the map $\sigma\colon C^\infty(\Omega) \to \Ge{\Omega}$ which assigns to $f \in C^\infty(\Omega)$ the class of 
\[ (\varphi, x) \mapsto f(x)\qquad ((\varphi,x) \in U(\Omega)). \]
Finally, partial derivatives extending those of distributions are given by
\[ (\pd_i R)(\varphi,x) \coleq \pd_i ( R ( \varphi,.)) (x). \]
Comparing this formula for $\pd_i R$ with \eqref{star} and the basic space $\Geb\Omega$ with $\Gob\Omega$, one sees that one does not need calculus on infinite-dimensional spaces anymore for the formulation of this algebra, hence the name \emph{elementary} algebra. Colombeau algebras possessing a canonical embedding of distributions, like $\SGo$ or $\SGe$, are traditionally called \emph{full} algebras.

\subsection{Special algebras}

A different simplification comes from fixing $\varphi$ in \eqref{genfunct} and viewing generalized functions as functions depending on $\e$ and $x$. This gives rise to the so-called \emph{simplified} (or \emph{special}) algebra, where generalized functions are represented by nets in $\e$ of smooth functions in $x$. This approach links particularly well with the sequential approach to distributions of Mikusi\'nski \cite{zbMATH03421400}. The definitions are as follows (\cite{GKOS, MOBook, 0601.35014}), where we employ the notation $K \csub \Omega$ for saying that $K$ is a compact subset of $\Omega$: 

\begin{definition}\label{def_specalg}
 We set
\begin{align*}
 I &\coleq (0,1], \\
 \Gsb\Omega &\coleq (C^\infty(\Omega))^I, \\
 \Gsm\Omega &\coleq \{ (u_\e)_\e \in \Gsb\Omega\ |\ \forall K \csub \Omega\ \forall \alpha \in \bN_0^n\\
& \qquad\qquad\qquad\qquad\qquad \exists N \in \bN: \sup_{x \in K}\abso{\pd^\alpha u_\e(x)} = O(\e^{-N}) \}, \\
 \Gsn\Omega &\coleq \{ (u_\e)_\e \in \Gsb\Omega\ |\ \forall K \csub \Omega\ \forall \alpha \in \bN_0^n\\
& \qquad\qquad\qquad\qquad\qquad \forall m \in \bN: \sup_{x \in K}\abso{\pd^\alpha u_\e(x)} = O(\e^m) \}, \\
  \Gs\Omega &\coleq \Gsm\Omega / \Gsn\Omega.
\end{align*}
\end{definition}

Smooth functions are trivially included as constant nets via
\begin{align*}
\sigma\colon C^\infty(\Omega) &\to C^\infty(\Omega)^I,\\
 f &\mapsto (f)_\e.
\end{align*}
In order to have $\iota$ and $\sigma$ agree on $C^\infty(\Omega)$, $\varphi$ has to satisfy all moment conditions defining the spaces $\cA_q(\bR^n)$ at once, which forces it to be an element of the Schwartz space $\cS(\bR^n)$. As a consequence, one can directly embed only compactly supported distributions at first and construct an embedding of $\cD'(\Omega)$ using the sheaf property afterwards (see \cite[Section 1.2.2]{GKOS} for details). However, one can also write this embedding in the form
\begin{align*}
 \iota u \coleq ( x \mapsto \langle u, \vec\psi_\e(x)\rangle)_\e \qquad (u \in \cD'(\Omega))
\end{align*}
for some net $(\vec\psi_\e)_\e \in C^\infty(\Omega, \cD(\Omega))^I$ (cf.~\cite{papernew}).

The special algebra is particularly convenient for direct extensions of operations on smooth functions because one can define them componentwise, i.e., for fixed $\e$. Spaces of generalized functions on manifolds and generalized sections of vector bundles are defined similarly to Definition \eqref{def_specalg}; for example, the space $\hat\cG^r_s(M)$ of generalized $(r,s)$-tensor fields on $M$ consists of nets $(t_\e)_\e \in (\cT^r_s(M))^I$ which are subject to the usual Colombeau-type quotient construction. A generalized metric then is an element of $\cG^0_2(M)$ such that each $t_\e$ is a pseudo-Riemannian metric of constant index (locally and for small $\e$) and which is non-degenerate in a certain sense. The development of the foundations of semi-Riemannian geometry in this setting came very naturally (\cite{genpseudo}). Analogously, one can define a space $\cG[\bR, M]$ of generalized mappings on $\bR$ with values in the manifold $M$ (\cite{Kgenmf}).

However, special algebras suffer from a major drawback: they fail to be diffeomorphism invariant and hence are only of limited use in a geometrical setting. Furthermore, the embedding $\iota$ does not commute with arbitrary Lie derivatives.

We will see in Section \ref{sec_funcanal} how these two essentially different approaches, full and special algebras, can be reconciled in a bigger picture.

\section{Applications to general relativity}\label{sec_appl}

The two papers \cite{genrel, JVdistgen} give a comprehensive overview of applications of distributions and generalized
functions to general relativity. Additionally to the works mentioned in Section \ref{sec_lindist}, these
articles discuss applications to Schwarzschild and Kerr space-times, ultrarelativistic black holes, conical singularities,
cosmic strings, weak singularities and generalized hyperbolicity. Therefore we only discuss more recent works, starting with
the description of an important class of wave-type space-times.

\subsection{Geodesic completeness of non-smooth space-times}
This class of space-times, the so-called \emph{$N$-fronted waves with parallel rays} (NPWs) for short, are defined
as a product $M=N\times \bR^2$, with metric
\begin{equation}\label{eq-int-met}
 l\ =\ \pi^*(h) + 2\ud u\ud v + F(.,u)\ud u^2\ ,
\end{equation}
where $h$ denotes the metric of an arbitrary, connected Riemannian manifold $(N,h)$, $\pi\colon M\rightarrow N$ is the projection
and $u,v$ are global null-coordinates on the $2$-dimensional Minkowski space $\bR^2_1$. Moreover, $F\in C^\infty(N\times\bR)$
is classically assumed to be smooth and is called the \emph{profile function} of this space-time. At first these models
were studied by Brinkmann in the context of conformal mappings of Einstein spaces (\cite{Brinkmann:25}). In the last ten
years NPWs were investigated with respect to their geometric properties and causal structure in \cite{CFS:03, FS:03, CFS:04,
FS:06}, where these space-times are called ``general plane waves''. Due to the geometric interpretation given in \cite{SS:12}
of $N$ as the wave surface of these waves, it seemed more natural to call them $N$-fronted waves, rather than plane-fronted
waves. Note that plane-fronted waves with parallel rays (pp-waves) (cf.~\cite[Chapter 17]{GrPo:09}) are a special case of NPWs
($N=\bR^2$ with the Euclidean metric).

Impulsive NPWs (iNPWs for short) are a generalization of impulsive pp-waves. Here the profile function
$F$ is given as a product $f\cdot\delta$, where $f\in C^\infty(N)$ and $\delta$ denotes the Dirac
distribution on the hypersurface $u=0$. As already shown in \cite{geodesics}, impulsive pp-waves can be thought of as being
geodesically complete in the sense that all geodesics (in the generalized sense) can be defined for all time. This
result also holds for the general class of iNPWs without any growth restrictions on the profile function $f$, which is in
contrast to the smooth case where on needs ``subquadratic'' behavior of the profile function in the $N$-component (see
\cite[Theorem 4.1]{FS:06}). Approximating the $\delta$ distribution by a \emph{strict delta net} $(\rho_\e)_\e$, we obtain a
family of smooth metrics
\[ l_\e = \pi^*(h) + 2 \ud u \ud v + f \rho_\e(u)\ud u^2. \]
Then the precise statement is the following (see \cite[Theorem 3.2]{SS:12}):

\begin{theorem}
Let $(N,h)$ be a connected and complete Riemannian manifold. Then for all $x_0\in N$, $\dot{x}_0\in T_{x_0}N$ and all $v_0$,
$\dot v_0\in{\bR}$ there exists $\e_0$ such that the solution $(x_\e,v_\e)$ of the geodesic equation with initial data
$x_\e(-1)=x_0$, $\dot x_\e(-1)=\dot x_0$, $v_\e(-1)=v_0$, $\dot v_\e(-1)=\dot v_0$ is defined for all $u\in\bR$, provided
$\e\leq\e_0$.
\end{theorem}

Note that $\e_0$ depends on the initial conditions and in general one cannot find an $\e_0$ that works for any initial
condition. One could say that ``in the limit as $\e\searrow 0$'' the space-time is complete. But without the use of Colombeau
generalized functions one cannot make this intuitive notion precise. Hence we arrive at the following definition (\cite{SS:13}):

\begin{definition}[Geodesic completeness for generalized metrics]
Let $g\in\cG^0_2(M)$ be a generalized metric of Lorentzian signature. Then the generalized space-time $(M,g)$ is said to be
\emph{geodesically complete} if every geodesic can be defined on $\bR$, i.e., every solution of the geodesic equation
\[ c''=0 \]
is in $\cG[\bR,M]$.
\end{definition}

In this sense iNPWs are geodesically complete (cf.~\cite{SS:13}) and one can see the need to use Colombeau generalized
functions to give a mathematically rigorous meaning to that notion. This is an example were one can observe that (just) using
distributions one is unable to formulate in a mathematical precise manner a physically interesting concept. The above
discussion provides a strong incentive to use Colombeau algebras when dealing with non-smooth space-times.

\subsection{Solution theory of the wave equation and generalized global hyperbolicity}
Another recent development is a solution theory for the Cauchy problem on non-smooth manifolds with \emph{weakly singular}
Lorentzian metrics (\cite{HKS:12}). In this work a generalization of the (smooth) metric splitting for globally
hyperbolic Lorentzian manifolds (\cite{BS:03, BS:05}) is employed to prove existence and uniqueness for the Cauchy problem
with compactly supported initial data in the special Colombeau algebra $\SGs$. This notion is called \emph{globally hyperbolic metric
splitting} and its definition is given in (\cite[Definition 6.1]{HKS:12}). The space-time $(M,g)$ splits in this case
as $M\cong\bR\times S$, where $S$ is a Cauchy hypersurface (for every representative of $g$). Then the Cauchy problem reads
\begin{align} \label{eq-cauchy}
\left\{\begin{aligned}
        \Box u &=0 \textrm{ on }M, \\
        u &= u_0 \textrm{ on }S, \\
        \nabla_{\xi} u &= u_1 \textrm{ on }S.
       \end{aligned}
       \right.
\end{align}
Here $u_0, u_1\in \Gs S$ have compact support and $\xi$ denotes the unit normal vector field of $S$.
\cite[Theorem 6.3]{HKS:12} states that \eqref{eq-cauchy} has a unique solution $u\in\Gs M$ if one assumes (additionally to $g$
being weakly singular) better asymptotics of the time-derivative of the metric as well as the spatial derivative of $g_{00}$,
i.e., condition (B) of \cite{HKS:12}. The inhomogeneous problem and the inclusion of lower-order terms can be handled
also by these methods (cf.~\cite{Hanel:11}).

This generalization of the smooth metric splitting is well-suited for dealing with the Cauchy problem but it is not directly
related to more geometric (existence of a Cauchy hypersurface) or more ``physical'' characterizations of global hyperbolicity
like satisfying the causality condition and having compact causal diamonds. Thus there was the need to see how well-suited
this notion of generalized global hyperbolicity is when dealing with concrete non-smooth space-times. In order to settle this
question, NPWs with non-smooth profile function were investigated in \cite{HS:13}. There the metric splitting was calculated
quite explicitly for a subclass of NPWs with smooth profile function and then these explicit constructions were used to
approximate the non-smooth metric splitting. One result obtained is that one should allow for a \emph{generalized
diffeomorphism} (in the sense of \cite[Section 4]{KS:99}) in the splitting of $M\cong \bR\times S$ (\cite[Proposition 3.3 and Corollary 3.5]{HS:13}). Another result is that if one approximates the non-smooth profile function appropriately then one obtains also an approximation of the metric splitting (\cite[Section 4]{HS:13}).

So the overall conclusion of \cite{HS:13} is that there are reasonable examples of non-smooth space-times possessing the
(adapted) globally hyperbolic metric splitting but still it would be desirable to have more
geometric generalizations of global hyperbolicity in the non-smooth setting.

In applications (especially to general relativity) it is often more convenient to use the special algebra, but of course it would be more reasonable to use diffeomorphism invariant Colombeau algebras, which are explained in the next section. It is our hope that the
functional analytic approach outlined in Section \ref{sec_funcanal} removes a lot of technicalities when dealing with these and hence enables one to take advantage of its merits more directly.

\section{Diffeomorphism invariance}\label{sec_diffinv}

Since the first introduction of Colombeau algebras there was the desire to have a functorial version of the theory. More explicitly, any diffeomorphism $\Omega \to \Omega'$ should induce an action $\cG(\Omega) \to \cG(\Omega')$ in a way respecting composition and the identity, and $\Omega$ should be allowed to be a manifold instead of an open subset of $\bR^n$. Such a diffeomorphism invariant construction would ensure independence of any choice of coordinate system and is essential for applications in a geometrical context.

\subsection{The diffeomorphism invariant local Colombeau algebra $\SGd$}

Accomplishing this goal was a long and tedious undertaking which spanned quite some years, involved several authors, and brought with it a lot of technical difficulties (\cite{colmani, Jelinek, found}; see \cite[Section 2.1]{GKOS} for more detailed comments).

The root of the problem lies in the fact that attempts focused on building a diffeomorphism invariant algebra on top of $\SGe$, but this algebra uses the linear structure of $\bR^n$ in an essential way for a simple and elementary presentation. As a consequence, it is particularly ill suited for a global formulation, which is why the necessary adaptions to obtain diffeomorphism invariance were so complicated. One had to 1) replace test objects $\varphi \in \cA_0(\bR^n)$ by ones depending also on $\e$ and $x$, 2) reintroduce calculus on infinite-dimensional locally convex spaces, 3) deal with poor domains of definition of test objects after being transformed by a diffeomorphism, 4) switch from vanishing moments to asymptotically vanishing moments and 5) establish invariance of moderateness and negligibility under derivatives, which became quite complicated in this setting.

Despite these difficulties, the diffeomorphism invariant local theory, substantially based on previous work of Colombeau and Meril \cite{colmani} and especially Jel\'inek \cite{Jelinek}, was for the first time obtained in \cite{found}, where the following definitions were given.

\begin{definition}\label{def_gd}We set $\Gdb\Omega \coleq C^\infty(U(\Omega))$.

Let $C^\infty_b(I \times \Omega, \cA_0(\bR^n))$ be the space of all smooth maps $\phi\colon I \times \Omega \to \cA_0(\bR^n)$ such that the corresponding map $\hat\phi\colon I \to C^\infty(\Omega, \cA_0(\bR^n))$ has bounded image.

$R \in \Gdb\Omega$ is called moderate if $\forall K \csub \Omega$ $\forall \alpha \in \bN_0^n$ $\exists N \in \bN$ $\forall \phi \in C^\infty_b(I \times \Omega, \cA_0(\bR^n))$: $\sup_{x \in K} \abso{\pd^\alpha ( R( S_\e \phi(\e,x), x))} = O(\e^{-N})$, where $\pd^\alpha ( R(S_\e\phi(\e,x),x))$ denotes the derivative of order $\alpha$ of the map $x' \mapsto R(S_\e \phi(\e,x'),x')$ at $x$. The set of all moderate elements will be denoted by $\Gdm\Omega$.

$R \in \Gdm\Omega$ is called negligible if $\forall K \csub \Omega$ $\forall m \in \bN$ $\exists q \in \bN$ $\forall \phi \in C^\infty_b(I \times \Omega,\cA_q(\bR^n))$: $\sup_{x \in K} \abso{R(S_\e\phi(\e,x),x)} = O(\e^m)$. The set of all negligible elements will be denoted by $\Gdn\Omega$.

We have embeddings $\iota\colon \cD'(\Omega) \to \Gdb\Omega$ and $\sigma\colon C^\infty(\Omega) \to \Gdb\Omega$ defined by
\begin{align*}
 (\iota u)(\varphi,x) &\coleq \langle u, \varphi(.-x) \rangle, \\
 (\sigma f)(\varphi, x) &\coleq f(x).
\end{align*}
Derivatives are given by
\[ (D_i R)(\varphi ,x)  \coleq \pd_i ( R (\varphi, .))(x). \]
The pullback of an element of $\Gdb{\Omega'}$ along a diffeomorphism $\mu\colon \Omega \to \Omega'$ is defined as $\mu^*R \in \Gdb\Omega$ given by
\[ (\mu^* R)(\varphi, x) \coleq R ( \mu_*(\varphi, x)), \]
where we set
\[ \mu_* (\varphi, x) \coleq (\varphi (\mu^{-1}(. + \mu x) - x) \cdot \abso{\det \D\mu^{-1}(.+\mu x)}, \mu x ). \]
Finally, $\Gd\Omega \coleq \Gdm\Omega / \Gdn\Omega$.
\end{definition}

A glance at the proofs of the basic constitutive properties of $\SGd$ (cf.~\cite{found}) shows that it is far from being simple and elegant. As will be seen below, it was in the global theory on manifolds where glimpses of a simpler formulation first became visible. On the base of this, the first author presented a short and concise treatment of $\SGd$ and $\SGh$ in \cite{gdnew}. But even with as much simplification of the proofs as possible, the problem remained that too many technicalities in its basic structure hindered the application and the further development of diffeomorphism invariant Colombeau algebras.

The main difficulties seem to lie in establishing diffeomorphism invariance of the test objects. We will see in Section \ref{sec_funcanal} how all of this can be replaced by a much more cleaner and simpler structure, but we will first sketch the further developments in the manifold setting based on $\SGd$.

\subsection{The global Colombeau algebra $\SGh$}

It is evident that the construction of $\SGd$ as given above does not directly generalize to the manifold setting because its formulation depends in an essential way on the linear structure of $\bR^n$, in particular concerning the test objects. One had to incorporate the formalism used by Jel\'inek in \cite{Jelinek} and find a coordinate-free characterization of the test objects in order to move to manifolds. In the following definition of the test objects for $\SGh$ (\cite[Definition 3.3]{global}), $\fX(M)$ denotes the space of smooth vector fields on a manifold $M$ and $B_r(p)$ the metric ball around $p \in M$ with radius $r$ with respect to an arbitrary (but fixed) metric.

\begin{definition}Set $\hat\cA_0(M) \coleq \{ \omega \in \Omega^n_c(M)\ |\ \int \omega = 1 \}$.

We denote by $\widetilde\cA_0(M)$ the set of all $\Phi \in C^\infty(I \times M, \hat\cA_0(M))$ satisfying
\begin{enumerate}[label=(\roman*)]
 \item $\forall K \csub M$ $\exists \e_0, C>0$ $\forall p \in K$ $\forall \e \le \e_0$: $\supp \Phi(\e,p) \subseteq B_{\e C}(p)$,
 \item $\forall K \csub M$ $\forall l,m \in \bN_0$ $\forall \zeta_1,\dotsc,\zeta_l,\theta_1,\dotsc,\theta_m \in \fX(M)$:
 \[ \sup_{p \in K, q \in X} \norm {L_{\theta_1} \dotsc L_{\theta_m}  ( L'_{\zeta_1} + L_{\zeta_1} ) \dotsc (L'_{\zeta_l} + L_{\zeta_l}) \Phi(\e,p)(q)} = O(\e^{-(n+m)}). \]
\end{enumerate}
where $L'_\zeta \Phi(\e,p)(q) \coleq L_\zeta ( p \mapsto \Phi(\e,p)(q))$. 

For $k \in \bN$ we denote by $\widetilde\cA_k(M)$ the set of all $\Phi \in \widetilde\cA_0(M)$ such that $\forall f \in C^\infty(M)$ and $\forall K \csub M$:
\begin{equation}\label{moment}
\sup_{p \in K} \abso{ f(p) - \int_M f(q) \Phi(\e,p)(q)} = O(\e^{k+1}).
\end{equation}
\end{definition}
Elements of $\widetilde\cA_0(M)$ were called ``smoothing kernels'' in \cite{global}, but we propose a different terminology where we will call these objects ``nets of smoothing kernels'' instead, the meaning of which will become clear in Section \ref{sec_funcanal}.

The spaces $\widetilde\cA_k(M)$ correspond, via scaling and translation (\cite[Lemma 4.2]{global}), to the spaces $C^\infty_b(I \times \Omega, \cA_q(\bR^n))$ of local test objects in $C^\infty_b$ from Definition \ref{def_gd}, but note that they still are not easily seen to be diffeomorphism invariant --- the proof of this property rests entirely on the local theory.

With this, one can define the space of scalar generalized functions on a manifold as follows:

\begin{definition}
We set $\Ghb M \coleq C^\infty(\hat\cA_0(M) \times M)$. $R \in \Ghb M$ is called moderate if $\forall K \csub M$ $\forall l \in \bN_0$ $\exists N \in \bN$ $\forall \zeta_1,\dotsc,\zeta_l \in \fX(M)$ $\forall \Phi \in \widetilde\cA_0(M)$: $\sup_{p \in K}\abso{L_{\zeta_1} \dotsc L_{\zeta_l}(R(\Phi(\e,p),p))} = O(\e^{-N})$, where $L_{\zeta_1} \dotsc L_{\zeta_l} (R(\Phi(\e,p),p))$ denotes the iterated Lie derivative of $p' \mapsto R(\Phi(\e, p'), p')$ at $p$. The subset of moderate elements of $\Ghb M$ is denoted by $\Ghm M$.

$R \in \Ghb M$ is called negligible if $\forall K \csub M$ $\forall l,m \in \bN_0$ $\exists k \in \bN$ $\forall \zeta_1,\dotsc,\zeta_l \in \fX(M)$ $\forall \Phi \in \widetilde\cA_k(M)$: $\sup_{p \in K} \abso{L_{\zeta_1} \dotsc L_{\zeta_l} ( R (\Phi(\e,p),p))} = O(\e^m)$. The subset of negligible elements of $\Ghb M$ is denoted by $\Ghn M$.

We have embeddings $\iota\colon \cD'(\Omega) \to \Ghb M$ and $\sigma\colon C^\infty(\Omega) \to \Ghb M$ defined by
\begin{align*}
 (\iota u)(\omega,p) &\coleq \langle u, \omega \rangle, \\
 (\sigma f)(\omega, x) &\coleq f(x).
\end{align*}
The Lie derivative of $R \in \Ghb M$ with respect to a vector field $X \in \fX(M)$ is defined as
\[ (\hat L_X R)(\omega, p ) \coleq -\ud_1 R(\omega, p)(L_X \omega) + L_X ( R (\omega, .))|_p, \]
where $\ud_1$ is the differential with respect to the first variable, and the action of a diffeomorphism $\mu \colon M \to N$ is given by
\[ (\mu_* R)(\omega,p) \coleq \mu_*( R(\mu^*\omega, \mu^{-1} p )). \]
We set $\Gh M \coleq \Ghm M / \Ghn M$.
\end{definition}

As in $\SGd$, pullback and Lie derivatives commute with the embeddings $\iota$ and $\sigma$ here. It should be noted that if one takes these definitions as the starting point also for $\SGd$ the local theory gets somewhat more transparent. An essential advantage of this approach then is that one has no problems at all with the domain of definition of generalized functions and of test objects, so part of the difficulties in $\SGd$ is simply avoided without any real drawbacks. Furthermore, the equivalence between local and global formalism becomes much more easier this way, as is evidenced in \cite{gdnew}. However, the difficulties with showing stability of moderateness and negligibility under diffeomorphisms and derivatives remains.

After the introduction of the global algebra $\SGh$ it was naturally very desirable to have an extension to a theory of generalized sections of vector bundles, in particular a theory of nonlinear generalized tensor fields for applications in (semi-)Riemannian geometry, but this faced some serious difficulties which we will describe next.

\subsection{Generalized tensor fields $\SGhrs$}

Coming from $\SGh$, the basic problem encountered in the tensor case is that one cannot use a coordinatewise embedding and simply take the definition $\Ghrs M \coleq \Gh M \otimes_{C^\infty(M)} \cT^r_s(M)$, where $\Gh M$ is the scalar algebra defined above. This cannot succeed due to a consequence of the Schwartz impossibility result \cite[Proposition 4.1]{global2}.

The underlying reason is that Colombeau algebras always involve some kind of regularization; this is visible either in the basic space and the embedding (special algebras) or in the testing procedure (full algebras), and in any case in the form embedded distributions take when they are being tested for moderateness or negligibility (see Section \ref{sec_funcanal} below). But in order to regularize non-smooth or distributional sections of a vector bundle one needs to transport vectors between different fibers of the bundle; this is intuitively clear if one tries to generalize to sections the usual method of regularizing by convolution with a smooth mollifier. On these grounds, the vector bundle of \emph{transport operators} was introduced in \cite[Appendix A]{global2} as follows:
\[ TO(M,M) \coleq \bigcup_{\mathclap{(p,q) \in M \times M}}\, \{ (p,q) \} \times L(T_pM, T_qN). \]
This means that the fiber over $(p,q) \in M \times M$ consists of all linear maps from $T_pM$ to $T_qN$, where $T_pM$ is the fiber over $p$ of the tangent bundle of $M$. Using a transport operator $A \in TO(M,M)$ and a smoothing kernel $\Phi \in C^\infty(M, \Omega^n_c(M))$ one can then approximate non-smooth vector fields $X$ by smooth ones $\tilde X$ as in
\[ \tilde X(p) \coleq \int_M A(q,p) X(q) \Phi(p)(q), \]
where, roughly speaking, $\tilde X$ will be close to $X$ if $A$ is close to the identity near the diagonal and $\Phi(p)$ is close to $\delta_p$ for all $p \in M$. This formula also applies to distributional $X$.

Based on these ideas, in \cite{global2} an algebra $\bigoplus_{r,s} \hat\cG^r_s$ of generalized tensor fields was constructed whose elements are represented by smooth functions $R(\omega, p, A)$ depending on a compactly supported $n$-form $\omega$, a point $p \in M$, and a transport operator $A \in TO(M,M)$.

This generalized tensor algebra, however, suffers from serious drawbacks. First of all, it inherits all the technical difficulties from $\cG^d$ and $\hat\cG$ and adds even more on top of it. Second, it is no sheaf: the corresponding proof which worked in all previous algebras breaks down due to the failure of the test objects to be `localizing' in a certain sense. And third, there was no way to define a covariant derivative $\nabla_X$ on $\SGhrs$ which is $C^\infty$-linear in the vector field $X$, which is an indispensable necessity for geometrical applications. The reason for this is that on scalars, the covariant derivative should agree with the Lie derivative, but the natural choice of the latter cannot be $C^\infty$-linear in $X$ again because of the Schwartz impossibility result. In sum, despite its achievements $\SGhrs$ still was unsatisfactory and raised the following questions:
\begin{enumerate}[label=(\roman*)]
 \item What is the deeper reason for $\SGhrs$ not to be a sheaf and not to allow for the introduction of a covariant derivative?
 \item Which adaptions are necessary for obtaining these features?
 \item Is there any way to obtain a more natural, less technical formulation of the whole theory?
\end{enumerate}

As we will see in the next section, there is indeed an approach which completely answers these questions in a very satisfactory way.

\section{The functional analytic approach}\label{sec_funcanal}

The shortcomings of $\SGhrs$, i.e., the missing sheaf property and the impossibility of devising a covariant derivative for it, brought the need for a more structured approach to Colombeau algebras, even more so because the previous ad-hoc approaches leading to $\SGd$, $\SGh$ and $\hat\cG^r_s$ gave no insight on how to proceed further.

In the following we will describe a functional analytic approach to Colombeau algebras (\cite{papernew}) which not only contains those algebras described above as special cases, but also gives a clear answer on how an algebra of generalized tensor fields with the desired properties can be successfully constructed. Furthermore, this setting seems to be very promising also for other aspects of Colombeau theory.

\subsection{Basic spaces}

The intuitive understanding of Colombeau algebras until now has been that some kind of smoothing has to be involved either in the embedding or in testing. Formally, this can be seen by examining the form embedded distributions take when they are being tested. The following table lists the corresponding expressions in the various Colombeau algebras we have considered so far.
\begin{center}
\begin{tabular}{cc}
 Variant & Test of $\iota u$, $u \in \cD'$\\
 \hline
 $\SGo$ & $\langle u, \varphi_{\e,x} \rangle$ \\
 $\SGe$ & $\langle u, (S_\e\varphi)(.-x) \rangle$ \\
 $\SGs$ & $\langle u, \vec\psi_\e(x)  \rangle$ \\
 $\SGd$ & $\langle u, S_\e\phi(\e,x)(.-x)  \rangle$ \\
 $\SGh$ & $\langle u, \Phi(\e,p)  \rangle$
\end{tabular}
\end{center}
Close inspection reveals that the argument of $u$ in every case is a test function depending on $\e$ and $x$ which converges to $\delta_x$ in a certain manner for $\e \to 0$. The precise link between this observation and the intuitive idea of \emph{smoothing} is given by a variant of Schwartz' kernel theorem (\cite[Th\'eor\`eme 3]{FDVV}), which reads
\[ L_b(\cD'(\Omega), C^\infty(\Omega)) \cong C^\infty(\Omega, \cD(\Omega)). \]
This means that if we want to regularize distributions in a reasonable (i.e., linear and continuous) way, this is done \emph{exactly} by letting them act on elements of $C^\infty(\Omega, \cD(\Omega))$, so-called \emph{smoothing kernels}. This space carries the topology of uniform convergence on compact sets in all derivatives. Explicitly, the correspondence between $\Phi \in L_b (\cD'(\Omega), C^\infty(\Omega))$ and $\vec\varphi \in C^\infty(\Omega, \cD(\Omega))$ is given by
\[
 \begin{aligned}
  (\Phi u)(x) & \coleq \langle u, \vec\varphi(x) \rangle\qquad &&(u \in \cD'(\Omega), x \in \Omega) \\
\vec\varphi(x) & \coleq \Phi^t (\delta_x) \qquad &&(x \in \Omega)
 \end{aligned}
\]

where $\Phi^t \in L_b(C^\infty(\Omega)', \cD(\Omega))$ is the transpose of $\Phi$. In light of this isomorphism, generalized functions in the sense of Colombeau in fact are, in the most encompassing sense, mappings from the space of smoothing kernels into the space of smooth functions. This viewpoint suggests the introduction of the basic space
\[ \Gfb\Omega \coleq C^\infty( \SK\Omega, C^\infty(\Omega)) \]
where $\SK\Omega \coleq C^\infty(\Omega, \cD(\Omega))$ is the space of smoothing kernels. This basic space trivially contains both distributions and smooth functions via the embeddings
\begin{equation}\label{embeddings}
\begin{aligned}
\iota\colon \cD'(\Omega) \to \Gfb\Omega&,\ \iota(u)(\vec\varphi)(x) \coleq \langle u, \vec\varphi(x) \rangle,\\
\sigma\colon C^\infty(\Omega) \to \Gfb\Omega&,\ \sigma(f)(\vec\varphi) \coleq f.
\end{aligned}
\end{equation}
We remark that this idea was also followed in \cite{1177.46033,1204.58021} on a more abstract level; what will be important for us is to connect this approach back to existing Colombeau theory.

It turns out very quickly that one needs to restrict this space in order to have any hope of obtaining sheaf properties. The simple observation that $(\iota u)(\vec\varphi)(x) = \langle u, \vec\varphi(x) \rangle$ depends only on $\vec\varphi(x)$ led to the introduction of the following so-called \emph{locality conditions}:

\begin{definition}\label{locdef}A function $R \in \Gfb\Omega = C^\infty(\SK\Omega, C^\infty(\Omega))$ is called
\begin{enumerate}[label=(\roman*)]
  \item \emph{local} if for all open subsets $U \subseteq \Omega$ and  smoothing kernels $\vec\varphi, \vec\psi \in \SK\Omega$ the equality $\vec\varphi|_U=\vec\psi|_U$ implies $R(\vec\varphi)|_U=R(\vec\psi)|_U$;
  \item \emph{point-local} if for all $x \in \Omega$ and smoothing kernels $\vec\varphi, \vec \psi \in \SK\Omega$ the equality $\vec\varphi(x)=\vec\psi(x)$ implies $R(\vec \varphi)(x)=R(\vec\psi)(x)$;
  \item \emph{point-independent} if for all $x,y \in \Omega$ and smoothing kernels $\vec\varphi, \vec\psi \in \SK\Omega$ the equality $\vec\varphi(x)=\vec\psi(y)$ implies $R(\vec \varphi)(x)=R(\vec \psi)(y)$.
\end{enumerate}
We denote by $\Gfbloc\Omega$, $\Gfbploc\Omega$ and $\Gfbpi\Omega$ the subsets of $\Gfb\Omega$ consisting of local, point-local and point-independent elements, respectively, and by $\Gflbloc\Omega$, $\Gflbploc\Omega$ and $\Gflbpi\Omega$ the corresponding subsets of $\Gflb\Omega \coleq L_b(\SK\Omega, C^\infty(\Omega))$.
\end{definition}

Now one can very easily see the following isomorphisms, which recover the basic spaces of $\cG^o$ and $\cG^d$ as well as distributions and distributions depending smoothly on a parameter as subspaces of $\cE(\Omega)$:
\begin{theorem}\label{subsp}
 \begin{enumerate}[label=(\roman*)]
  \item $\Gfbpi\Omega \cong C^\infty(\cD(\Omega))$,
  \item $\Gfbploc\Omega \cong C^\infty(\cD(\Omega), C^\infty(\Omega))$,
  \item $\Gflbpi\Omega \cong \cD'(\Omega)$,
  \item $\Gflbploc\Omega \cong C^\infty(\Omega, \cD'(\Omega))$.
 \end{enumerate}
\end{theorem}

The importance of Theorem \ref{subsp} lies in the fact that $\cE(\Omega)$ is not simply one more of many possible basic spaces, but in a sense the most general one, and still contains common basic spaces as subalgebras. We will see below that using $\SGfbloc$ (or any subalgebra of it) as a basic space, the resulting Colombeau algebra will be a sheaf if the right test objects are used for testing. This was not the case for $\SGhrs$, which is why it fails to be a sheaf.

\subsection{Testing}

Having identified $\Gfb\Omega$ as the basic space underlying all Colombeau algebras, the next question was how the quotient construction could be transferred to it in a natural way. Suppose we are given a net of smoothing kernels $(\vec\varphi_\e)_\e \in \SK\Omega^I$, corresponding to a net of smoothing operators $(\Phi_\e)_\e \in L_b(\cD'(\Omega), C^\infty(\Omega))^I$. What conditions do we need for the usual Colombeau-type quotient construction?

The intuitive understanding of test objects always has been that asymptotically, they behave like delta functionals. This is seen clearly in the local prototypical smoothing kernels $\varphi_{\e,x}(y) \coleq \e^{-n} \varphi ( (y-x) / \e )$ with $\varphi \in \cA_q(\bR^n)$. Using the language of smoothing operators, this asymptotic behavior can be formulated precisely as follows:
\begin{equation}\label{sk_1}
\Phi_\e \to \id \textrm{ in }L_b(\cD'(\Omega), \cD'(\Omega)).
\end{equation}
The moment condition \eqref{locmoment}, or \eqref{moment} on a manifold, corresponds to convergence to the identity $\id_{C^\infty(\Omega)}$ on $C^\infty(\Omega)$ in the following sense:
\begin{equation}\label{sk_2}
\forall p \in \cs(L_b(C^\infty(\Omega), C^\infty(\Omega)))\ \forall m \in \bN: p (\Phi_\e|_{C^\infty} - \id ) = O(\e^m).
\end{equation}
Here, $\cs(E)$ denotes the set of continuous seminorms of a locally convex space $E$. Note that here we abolished the grading on the space of test objects known from full algebras, which is a relict of the fact that a nonzero test function $\varphi \in \cD(\Omega)$ cannot have all its moments vanishing at once. Instead, we use the kind of convergence used in special algebras with our full basic space.

These conditions do not yet guarantee that embedded distributions are moderate; for this, we furthermore have to impose that
\begin{equation}\label{sk_3}
\forall p \in \cs ( L_b(\cD'(\Omega), C^\infty(\Omega)))\ \exists N \in \bN: p (\Phi_\e) = O(\e^{-N}).
\end{equation}
\eqref{sk_1}, \eqref{sk_2} and \eqref{sk_3} are the essential conditions test objects have to satisfy for the quotient construction to work, in general. We shall add another (optional) condition which ensures the sheaf property if we test on a basic space whose elements are at least local (cf.~\cite{papernew}):
\begin{definition}\label{def_localnet} A net $(\vec\varphi_\e)_\e \in \SK\Omega^I$ is called \emph{localizing} if $\forall x \in \Omega$ $\exists$ an open neighborhood $V$ of $x$ $\forall r>0$ $\exists \e_0 \in \bN$ $\forall \e \le \e_0$ $\forall x \in V$: $\supp \vec\varphi_\e(x) \subseteq \overline{B_r}(x)$.
\end{definition}

Using these ingredients, one can build a diffeomorphism invariant Colombeau algebra without much effort, as we will see in the next section.

Concerning the introduction of a covariant deriative of generalized tensor fields, we have mentioned that the problem essentially lies in obtaining a $C^\infty$-linear Lie derivative. Because the natural Lie derivative (which is obtained by differentiating with respect to time the pullback along the flow of a given vector field) commutes with the embedding $\iota$ it cannot be $C^\infty$-linear due to the Schwartz impossibility result, hence one needs to introduce another Lie derivative $\widetilde L_X$ which is $C^\infty$-linear and somehow related to the natural one. The solution is to set $(\widetilde L_X R)(\vec\varphi) \coleq L_X ( R (\vec\varphi))$. In a sense, this means that one differentiates after regularization. This is exactly the approach which is used in special Colombeau algebras, which is why the definition of a covariant derivative was without problems there (\cite{genpseudo}). While in the case of $\SGh$ and $\SGhrs$ the map $\widetilde L_X$ cannot be defined (it maps out of the space of point-local elements), this is made possible in our approach by the bigger basic space. Finally, it can be seen that the natural Lie derivative and $\widetilde L_X$ agree for embedded distributions \emph{on the level of association}, which justifies this definition.

\subsection{A new diffeomorphism invariant algebra of generalized functions}

We gave a treatment of both the basic space and the testing procedure in functional analytic terms. In this section we will use this to propose a replacement for $\cG^d$ of \cite{found} which is considerably simpler and directly translates to manifolds and, with some modifications, also to the vector valued case.

We will base our presentation on $\Gfb\Omega$ in order to show the general scheme, but one can do the same without problems on $\Gfbloc\Omega$ in order to obtain a sheaf, or even on $\Gfbploc\Omega$ or $\Gfbpi\Omega$ as long as one only performs operations preserving (point-)locality or point-independence, respectively.

Hence, the basic definitions take the following form.

\begin{definition}
 We define the basic space $\Gnb\Omega$ and embeddings $\iota\colon \cD'(\Omega) \to \Gnb\Omega$, $\sigma\colon C^\infty(\Omega) \to \Gnb\Omega$ by
\begin{align*}
 \SK\Omega & \coleq C^\infty(\Omega, \cD(\Omega)), \\
 \Gnb\Omega & \coleq C^\infty ( \SK\Omega, C^\infty(\Omega)), \\
 (\iota u)(\vec\varphi)(x) &\coleq \makebox[10ex][l]{$\langle u, \vec\varphi(x) \rangle$} (u \in \cD'(\Omega), \vec\varphi \in \SK\Omega, x \in \Omega),\\
 (\sigma f)(\vec\varphi) &\coleq \makebox[10ex][l]{$f$} (f \in C^\infty(\Omega), \vec\varphi \in \SK\Omega).
\end{align*}
Given a diffeomorphism $\mu \colon \Omega \to \Omega'$, its action $\mu_*\colon \Gnb\Omega \to \Gnb{\Omega'}$ is defined as
\[
(\mu_* R)(\vec\varphi) \coleq \mu_* ( R ( \mu^*\vec\varphi)) = R(\mu^*\vec\varphi) \circ \mu^{-1},
\]
where $(\mu^* \vec\varphi)(x) \coleq (\vec\varphi(\mu x) \circ \mu) \cdot \abso{\det \D\mu}$ is the natural pullback of smoothing kernels. The directional derivative with respect to a vector field $X \in C^\infty(\Omega, \bR^n)$ is defined as
\[  (L_X R) (\vec\varphi) \coleq  - (\ud R) (\vec\varphi) (L^{SK}_X \vec\varphi) + L_X ( R ( \vec\varphi)), \]
where $L^{SK}_X\vec\varphi = L_X \circ \vec\varphi + L_X\vec\varphi $ is the Lie derivative of smoothing kernels.
\end{definition}
It follows directly from the definitions that $\iota$ and $\sigma$ commute with diffeomorphisms. For testing, we need the following spaces of \emph{test objects}:
\begin{definition}Let $S$ be the set of all nets $(\vec\varphi_\e)_\e \in C^\infty(\Omega, \cD(\Omega))^I$ such that the corresponding sequence $(\Phi_\e)_\e \in L_b(\cD'(\Omega), C^\infty(\Omega))^I$ satisfies, for $\e \to 0$,
 \begin{enumerate}[label=(\roman*)]
  \item $\Phi_\e \to \id \textrm{ in }L_b(\cD'(\Omega), \cD'(\Omega))$,
  \item $\forall p \in \cs(L_b(C^\infty(\Omega), C^\infty(\Omega)))\ \forall m \in \bN: p (\Phi_\e|_{C^\infty(\Omega)} - \id_{C^\infty(\Omega)} ) = O(\e^m)$,
  \item $\forall p \in \cs ( L_b(\cD'(\Omega), C^\infty(\Omega)))\ \exists N \in \bN: p (\Phi_\e) = O(\e^{-N})$.
 \end{enumerate}
Let $S_0$ be the set of all $(\vec\varphi_\e)_\e$ such that $(\Phi_\e)_\e$ satisfies
\begin{enumerate}[label=(\roman*)]
 \item $\Phi_\e \to 0 \textrm{ in }L_b(\cD'(\Omega), \cD'(\Omega))$,
  \item $\forall p \in \cs(L_b(C^\infty(\Omega), C^\infty(\Omega)))\ \forall m: p (\Phi_\e|_{C^\infty} ) = O(\e^m)$,
  \item $\forall p \in \cs ( L_b(\cD'(\Omega), C^\infty(\Omega)))\ \exists N \in \bN: p (\Phi_\e) = O(\e^{-N})$.
\end{enumerate}
We may demand in both cases that $(\Phi_\e)_\e$ is localizing as in Definition \ref{def_localnet}.
\end{definition}
Note that for $(\vec\varphi_\e)_\e \in S$, $(L^{SK}_X\vec\varphi_\e)_\e$ will be an element of $S_0$. 
\begin{definition}An element $R \in \Gnb\Omega$ is called \emph{moderate} if $\forall p \in \cs (C^\infty(\Omega))$ $\forall k \in \bN_0$ $\exists N \in \bN$ $\forall (\vec\varphi_\e)_\e \in S$, $(\vec\psi^1_\e)_\e \dotsc (\vec\psi^k_\e)_\e \in S_0$:
\[ p \left(\, (\ud^k R) (\vec\varphi_\e)(\vec\psi_\e^1,\dotsc,\vec\psi_\e^k)\, \right) = O(\e^{-N}). \]
 The set of all moderate elements of $\Gnb\Omega$ is denoted by $\Gnm\Omega$.
 
 An element $R \in \Gnb\Omega$ is called \emph{negligible} if $\forall p \in \cs (C^\infty(\Omega))$ $\forall k \in \bN_0$ $\forall m \in \bN$ $\forall (\vec\varphi_\e)_\e \in S$, $(\vec\psi^1_\e)_\e \dotsc (\vec\psi^k_\e)_\e \in S_0$:
\[ p \left(\, (\ud^k R) (\vec\varphi_\e)(\vec\psi_\e^1,\dotsc,\vec\psi_\e^k)\, \right) = O(\e^m). \]
The set of all negligible elements of $\Gnb\Omega$ is denoted by $\Gnn\Omega$.
 
We set $\Gn\Omega \coleq \Gnm\Omega / \Gnn\Omega$.
\end{definition}
We summarize the main properties of these definitions in the following theorem:
\begin{theorem}
 \begin{enumerate}[label=(\roman*)]
  \item $\iota$ and $\sigma$ map into $\Gnm\Omega$;
  \item $\iota|_{C^\infty(\Omega)} - \sigma$ maps into $\Gnn\Omega$;
  \item $\iota$ is injective into $\Gn\Omega$;
  \item $\iota$ and $\sigma$ commute with diffeomorphisms and Lie derivatives;
  \item sums, products, diffeomorphisms and Lie derivatives preserve moderateness and negligibility and hence are well-defined on $\Gn\Omega$;
  \item the sets $S$ and $S_0$ of test objects are non-empty.
 \end{enumerate}
\end{theorem}

These properties follow directly from the definitions without any effort. We finish with some remarks:

\begin{enumerate}[label=(\roman*)]
 \item If we replace the basic space with $\SGfbloc$, $\SGfbploc$ or $\SGfbpi$ and test only with localizing elements of $S$ and $S_0$ then $\SGn$ is a sheaf.
 \item For the vector-valued case, one essentially only has to replace smoothing operators $L_b(\cD'(\Omega), C^\infty(\Omega))$ by the corresponding vector smoothing operators. This will be the subject of an upcoming paper of the first author.
 \item Instead of taking all test objects $S$ (and $S_0$) as above, one can restrict to certain classes of test objects. The properties of the algebra will directly depend on the choice of the class of test objects. Using only one fixed test object, for example, would give the classical special algebra. 
\end{enumerate}

\section{Conclusion}

We have seen that important physical situations lead to problems involving multiplication of distributions. While linear distribution theory in general fails in this setting, the theory of algebras of generalized functions in the sense of Colombeau leads to new, physically meaningful results. However, in order to obtain truly geometric results one needs to work in diffeomorphism invariant algebras, which were for a long time considered to be too technical and not well suited for applications. The approach we presented here is hoped to alleviate this shortcoming and to provide a better basis for further applications of Colombeau theory in general relativity and geometry, and also to provide a functional analytic foundation for and a better understanding of the field of Colombeau algebras.

{\bf Acknowledgments.} This work was supported by grants P23714 and P25326 of the Austrian Science Fund (FWF) and OEAD WTZ Project CZ15/2013.

\end{document}